# The bilinear maximal functions map into
## $L^p$ **for** $2/3 < p \le 1$

By Michael T. Lacey*

### Abstract

The bilinear maximal operator defined below maps $L^p \times L^q$ into $L^r$ provided $1 < p, q < \infty$, $1/p + 1/q = 1/r$ and $2/3 < r \le 1$.

$$Mfg(x) = \sup_{t>0} \frac{1}{2t} \int_{-t}^{t} |f(x+y)g(x-y)| \ dy.$$

In particular $Mfg$ is integrable if $f$ and $g$ are square integrable, answering a conjecture posed by Alberto Calderón.

## 1. Principal results

In 1964 Alberto Calderón defined a family of maximal operators by

$$Mfg(x) = \sup_{t>0} \frac{1}{2t} \int_{-t}^{t} |f(x-\alpha y)g(x-y)| \ dy, \qquad \alpha \neq 0, 1$$

which have come to be known as bisublinear maximal functions. He raised the striking conjecture that $Mfg$ is integrable if $f$ and $g$ are square integrable. A proof of this and more is provided in this paper.

1.1. Theorem. *Let* $\alpha \neq 0, 1$, *and let* $1 < p, q < \infty$ *and set* $1/r = 1/p + 1/q$. *If* $2/3 < r \le 1$ *then* $M$ *extends to a bounded map from* $L^p \times L^q$ *into* $L^r$.

Now, if $r > 1$, $M$ maps into $L^r$, as follows from an application of Hölder's inequality in the $y$ variable. Thus the interest is in the case $2/3 < r \le 1$. That $r$ can be less than one is intriguing and unexpected.

---

*This work has been supported by an NSF grant, DMS-9706884.

Our proof forsakes the maximal function for the maximal truncations of singular integrals. Let $K(y)$ be a singular integral kernel satisfying

$$(1.2) \qquad |K(y)| \quad \leq \quad \frac{C}{|y|},$$

$$(1.3) \qquad |\partial^n K(y)| \quad \leq \quad \frac{C}{|y|^{n+1}}, \quad n = 0, 1, \dots, N$$

where $N$ is some large integer, and

$$(1.4) \qquad \hat{K}(\xi) = \int e^{-2\pi i \xi y} K(y) \, dy \text{ is a bounded function of } \xi.$$

Kernels $K$ define bilinear operators by

$$Tfg(x) = \int f(x - \alpha y)g(x - y)K(y) \, dy, \qquad \alpha \neq 0, 1.$$

For the kernels of interest to us, the integral is defined *a priori* only for, say, functions $f$ and $g$ in the Schwartz class. Yet, the methods of [4], [5] prove that $T$ is bounded from $L^p \times L^q$ into $L^r$, provided $2/3 < r < \infty$. We extend this result for the maximal truncations as follows.

1.5. THEOREM. *There is an integer $N$ such that the following holds. Let $\alpha \neq 0, 1$, $1 < p, q \leq \infty$, and $1/r = 1/p + 1/q$. If $2/3 < r < \infty$, then for all Calderón-Zygmund kernels $K$ satisfying* (1.2)–(1.4), *the maximal operator below extends to a bounded map from $L^p \times L^q$ into $L^r$:*

$$T^* fg(x) = \sup_{\epsilon < \delta} \left| \int_{\epsilon < |y| < \delta} f(x - \alpha y)g(x - y)K(y) \, dy \right|.$$

The relationship between these two theorems is as follows. Let $\mathbf{T}^\bullet$ be defined as $\mathbf{T}^*$ is above, but with $\epsilon$ and $\delta$ restricted to the values $\{2^k \mid k \in \mathbf{Z}\}$. Then, we can choose two kernels $K_1$ and $K_2$ satisfying the hypotheses (1.2)–(1.4), so that for some constant $c$,

$$cMfg(x) \leq \int_{-1}^1 |f(x - \alpha y)g(x - y)| \, dy + T_1^\bullet fg(x) + T_2^\bullet fg(x),$$

where $\mathbf{T}_j^\bullet$ is specified by $K_j$. Our proof shows that $\mathbf{T}^\bullet$ satisfies the desired norm inequalities. Therefore, so does $M$. And conversely, we have $\mathbf{T}^* fg < c' Mfg + \mathbf{T}^\bullet fg$.

The method of proof is quite close to that of the author's prior collaboration with C. Thiele, [4, 5] yet the differences manifest themselves in two ways. The discussion in the prior work is done at a quite general level, namely the methods therein work equally well for the bilinear Hilbert transform and natural analogs in the Walsh-Paley setting. (See [9].) The nature of the maximal function forces us to abandon that level of generality. And more interestingly, the essence of the matter in [4] lies in the formulation and proof of certain almost orthogonal results. Maximal forms of these results must be proved; indeed this is essentially the only matter that we need address in this paper. These maximal inequalities rely in an essential way on a novel maximal inequality proved by Bourgain, [3]. Bourgain's inequality also plays a critical role in [2], a paper which has already demonstrated the close connection between the bilinear maximal inequality and Bourgain's lemma.

In the next section we define a class of operators, the "model sums," which are the principal concern of the proofs in [4], [5]. The bound for the model sums is stated. One may ask if the bound $r > 2/3$ is necessary in our theorems above. At least for the model sums, we show by example that our methods of proof do not apply to this case.

As is shown in [6], the bilinear singular integrals can be written as a sum of model sums. But the argument requires a degree of smoothness in the Calderón-Zygmund kernel $K$. This is the source of the extra derivatives imposed in Theorem 1.5, through hypothesis (1.4). The argument in [6] will go through if $N = 100$, a number which is certainly not sharp.

And then we turn to the maximal variants of the orthogonality results of [4], the subject which takes up most of this paper.

## 2. Model sums

The analysis of the operations in question requires overlapping representations of the functions involved, one representation for each scale of the problem, with the interaction between scales in representations manifesting itself in combinatorial ways. We define certain sums which facilitate a precise discussion of these issues. Unfortunately, the description of these sums is involved and phrased in general terms, with the generality forced upon us by the largely technical requirement of recovering the bilinear singular integrals. Other aspects of our definitions are forced upon us by the particular form of the maximal inequality (Lemma 3.9) below. The most concrete form of the Definition 2.8 below appears in [7] and all of the essential mathematical issues already appear in them.

2.1. *Definitions and theorem.* The combinatorics that enters into the problem is that of the dyadic rectangles in the space-frequency plane. Because of this connection, we will identify certain functions in the Schwartz class with rectangles in that plane, with the geometry of the rectangles encoding orthogonality properties of the functions. That is the purpose of these definitions.

Take the Fourier transform to be

$$\hat{f}(\xi) = \mathcal{F}f(\xi) = \int_{-\infty}^{\infty} e(-x\xi)f(x)\,dx$$

where $e(x) = e^{2\pi i x \xi}$. We also adopt the notation $\langle f, g \rangle = \int f\bar{g}\,dx$.

2.1. *Definition.* Let $\mathcal{R}$ be a class of rectangles, $\Phi = \{\varphi_\rho \mid \rho \in \mathcal{R}\}$ a set of Schwartz functions and $C_m, m \geq 1$ constants. $\Phi$ is said to be *adapted to $\mathcal{R}$ with constants $C_m$* if for all $\rho = I_\rho \times \omega_\rho \in \mathcal{R}$ the following conditions hold.

$$(2.2) \qquad\qquad\qquad \|\varphi_\rho\|_2 = 1,$$

$$(2.3) \qquad\quad |\rho| = |I_\rho| \times |\omega_\rho| \leq \sqrt{2}\,|\rho'| \ \text{ for all } \rho' \in \mathcal{R},$$

there is an affine function $a(\xi)$ so that

$$(2.4) \qquad\quad \widehat{\varphi_\rho}(\xi) \text{ is supported on } \frac{3}{4}\{a(\xi) \mid \xi \in \omega_\rho\},$$

$$(2.5) \qquad |\varphi_\rho(x)| \leq C_m \frac{1}{\sqrt{|I_\rho|}}\left(1 + \frac{|x - c(I_\rho)|}{|I_\rho|}\right)^{-m} \quad m \geq 0, x \in \mathbf{R},$$

where $c(J)$ denotes the center of the interval $J$, and

$$(2.6) \qquad\quad \text{if } \omega_\rho = \omega_{\rho'} \text{ but } I_\rho \neq I_{\rho'} \text{ then } \langle \varphi_\rho, \varphi_{\rho'} \rangle = 0.$$

The most obvious way to define the $\varphi_\rho$ is to begin with a fixed Schwartz function $\varphi$ of compact support in frequency and change its location, scale and frequency modulation according to the location of $\rho$. Thus,

$$\varphi_\rho(x) = e(c(\omega_\rho)x)\frac{1}{\sqrt{|I_\rho|}}\varphi\left(\frac{x - c(I_\rho)}{|I_\rho|}\right).$$

This is in fact what is done to pass from the model sums defined below to the bilinear singular integrals.

The first condition in the definition is just a normalization; the second condition arises from the necessary condition that $I_\rho$ and $\omega_\rho$ be in Fourier duality; and the next three conditions show that $\rho$ describes the location of $\varphi_\rho$ in the space-frequency plane. The requirements of this definition are more stringent than those of Definition 2.1 of [4].

We define the sums which model the bilinear operators and as well facilitate our combinatorial analysis.

2.7. *Definition.* Call a collection $\mathbf{I}$ of intervals in $\mathbf{R}$ a *grid* if for all $I \in \mathbf{I}$ we have $2^k \leq |I| \leq \frac{4}{3}2^k$ for some integer $k$, and if for all $I, I' \in \mathcal{I}$ we have $I \cap I' \in \{\emptyset, I, I'\}$.

2.8. *Definition.* A *model sum* is built up from three collections of rectangles $\mathcal{R}_i \subset \{I \times \omega i \mid I \in \mathcal{I}, \omega i \in \Omega_i\}$ for $i = 1, 2, 3$. These three collections are indexed by the same set $\mathbf{S}$, and for each $s \in \mathbf{S}$, the associated tiles in $\mathcal{R}_i$ all have the same first coordinate. Thus we write $\mathcal{R}_i = \{I_s \times \omega i(s) \mid s \in \mathbf{S}\}$. Let $\{\phi_{s,i} = \phi_{I_s \times \omega i(s)} \mid s \in \mathbf{S}\}$ be a set of Schwartz functions adapted to $\mathcal{R}_i$ with constants $C_m$. Assume that

$$(2.9) \qquad\qquad \mathcal{I} = \{I_s \mid s \in \mathbf{S}\} \text{ is a grid,}$$

$$(2.10) \qquad\qquad \Omega = \Omega_1 \cup \Omega_2 \cup \Omega_3 \text{ is a grid,}$$

$$(2.11) \qquad \sup\{\sigma \mid \xi \in \omega i(s)\} > \sup\{\sigma \mid \xi \in \omega j(s)\} \text{ for all } 1 \le i < j \le 3,$$

$$(2.12) \quad \text{if } \omega i(s) \subset_{\neq} \omega \text{ for any } s \in \mathbf{S} \text{ and } \omega \in \Omega, \text{ then } \omega j(s) \subset \omega \text{ for all } j.$$

$$(2.13) \qquad\qquad \operatorname{dist}(\omega i(s), \omega j(s)) \le K \left| \omega i(s) \right|, \quad s \in \mathbf{S}, i \neq j.$$

The model sum is then

$$\mathcal{M}(f_1, f_2)(x) = \sum_{s \in \mathbf{S}} \frac{\varepsilon_s}{\sqrt{|I_s|}} \prod_{i=1}^{2} \langle f_i, \phi_{s,i} \rangle \phi_{s,3}(x),$$

where $\varepsilon_s$ is an arbitrary choice of sign, $\varepsilon_s \in \{\pm 1\}$. The maximal form of these sums is

$$\mathcal{M}^{\max}(f_1, f_2)(x) = \sup_{k \in \mathbf{Z}} \Big| \sum_{\substack{s \in \mathbf{S} \\ |I_s| \ge 2^k}} \frac{\varepsilon_s}{\sqrt{|I_s|}} \prod_{i=1}^{2} \langle f_i, \phi_{s,i} \rangle \phi_{s,3}(x) \Big|.$$

The principal result of the paper asserts that model sums map into $L^r$, provided $2/3 < r < \infty$ while this result fails if $r < 2/3$.

2.14. **Theorem.** *Let $1 < p, q \le \infty$ and set $r = 1/p + 1/q$. Let $\mathcal{M}$ be a model sum.*

- *If $2/3 < r < \infty$, then $\mathcal{M}^{\max}$ maps $L^p \times L^q$ into $L^r$. The norm of $\mathcal{M}$ depends upon the choice of affine map in (2.4) and the constants in (2.5) but is otherwise independent of the choice of model sum.*

- *But for $r < 2/3$ there is a model sum that does not extend to a bounded bilinear map from $L^p \times L^q$ into $L^r$.*

The methods employed here do not give much of a clue as to what happens in the case of $p$ or $q$ being 1, and they would have to substantially refined to decide the case of $r = 2/3$.

We make a simple observation. In the definition of a model sum, let $\omega_s$ denote the convex hull of the three intervals $\omega_i(s)$ for $i = 1, 2, 3$. Then one readily checks that $\{\omega_s, \omega_i(s) \mid s \in \mathbf{S}, i = 1, 2, 3\}$ is a union of $O(1)$ grids. We can therefore identify the index set $\mathbf{S}$ with the rectangles $I_s \times \omega_s$. We write $si = I_s \times \omega i(s) = I_s \times \omega_{si}$.

The symmetry in the definition of the model sums is a central aspect of the proofs in [4]. The maximal function breaks this symmetry, of course, but it can be regained by the introduction of stopping times. Take $\sigma(x) : \mathbf{R} \to \{2^k \mid k \in \mathbf{Z}\}$ to be an arbitrary measurable function and set

$$\mathcal{M}^\sigma(f_1, f_2)(x) = \sum_{s \in \mathbf{S}} \frac{\varepsilon_s}{\sqrt{|I_s|}} \prod_{i=1}^{2} \langle f_i, \phi_{s,i} \rangle \phi_{s,3}^\sigma(x)$$

where

$$\phi_{s,3}^\sigma(x) = \phi_s(x) \mathbf{1}_{\{\sigma(x) \le |I_s|\}}.$$

Recall from Definition 2.7 that $1 \le |I_s| 2^k \le \frac{4}{3}$ for some integer $k$, so that it suffices to control $\mathcal{M}^\sigma$. With the inequalities of the next section, one can then analyze $\mathcal{M}^\sigma$ just as in [4], [5], thereby proving the theorem.

2.2.*Counterexample.* We prove the negative half of Theorem 2.14 by assuming that the indices $p$ and $q$ are such that the the model sums are unconditionally bounded as bilinear maps from $L^p \times L^q$ into $L^r$, then showing that necessarily $r \ge 2/3$.

An example model sum accomplishes this for us. The example is straightforward, being already apparent at a single scale. We take the rectangles associated to the model sums to be $\rho_{n,i} = [0,1) \times [4n+i, 4n+i+1)$ for $n \in \mathbf{Z}$ and $i = 1, 2, 3$. To define the functions adapted to these rectangles let $\phi$ be a Schwartz function with $\|\phi\|_2 = 1$ and $\hat{\phi}(\xi)$ supported on $[-1/4, 1/4]$. Then set

$$\phi_{n,i}(x) = e((4n+i+1/2)x)\phi(x).$$

Clearly $\phi_{n,i}$ is adapted to $\rho_{n,i}$. Moreover we have the model sums

$$\mathcal{M}_N(f_1, f_2)(x) = \sum_{n=0}^{N-1} \varepsilon_n \prod_{i=1}^{2} \langle f_i, \phi_{n,i} \rangle \phi_{n,3}(x), \quad \varepsilon_n \in \{\pm 1\}.$$

The numerous conditions which a model sum must satisfy are trivial to check in this instance.

Assume that the $\mathcal{M}_N$ are uniformly bounded, over all choices of sign and all $N$, as a map into $L^r$. Take for $i = 1, 2$, the functions $f_i$ to be

$$\begin{aligned}
f_i(x) &= \sum_{n=0}^{N-1} \phi_{n,i}(x) \\
&= e((i+1/2)x)\phi(x) \sum_{n=0}^{N-1} e(nx) \\
&= e((i+1/2)x)\phi(x) \frac{e(Nx)-1}{e(x)-1}.
\end{aligned}$$

Then we have the inequality

$$\|\mathcal{M}_N(f_1, f_2)\|_r = \left\|\sum_{n=0}^{N-1} \varepsilon_n \phi_{n,3}(x)\right\|_r \le K \|f_1\|_p \|f_2\|_q.$$

We average the left-hand side over all choices of signs and apply Khintchine's inequality. From this it follows that

$$\sqrt{N} \le K \|f_1\|_p \|f_2\|_q.$$

Yet is is easy to see that $\|f_i\|_t \le K' N^{1-1/t}$ for $1 < t < \infty$, so that we see the inequality $\sqrt{N} \le K N^{2-1/p-1/q} = K N^{2-1/r}$ valid for all $N$. Hence $r \ge 2/3$ is necessary for the uniform boundedness of the model sums.

## 3. Almost-orthogonality

We describe the orthogonality result upon which the theory of bilinear singular integrals is based. Then the maximal variant is stated, with the bulk of this section being taken up with its proof. The proof is quite combinatorial, with the principal aim being to arrange the functions involved into sets to which we can apply the critical maximal inequality of Bourgain, Lemma 3.9 below.

3.1. *Definition.* For two rectangles $I \times \omega$ and $I' \times \omega'$ we write $I \times \omega < I' \times \omega'$ if $I \subset I'$ and $\omega' \subset \omega$. A set of tiles $\mathbf{T}$ is call a 1-*tree with top* $t$ if $T = \{t\}$ or $s < t$ but $s1 \cap t = \emptyset$ for all $s \in \mathbf{T}$.

Let us set notation. Collections of tiles $\mathbf{S}$ are presumed to be unions of 1-trees $\mathbf{T}_t$ with tops $t \in \mathbf{S}^*$. The collection $\mathbf{S}^*$ need not be a subset of $\mathbf{S}$. Set

$$N_{\mathbf{S}}(x) = \sum_{t \in \mathbf{S}^*} \mathbf{1}_{I_t}(x).$$

Introduce a "stopping time" $\sigma : \mathbf{R} \to \{2^k \mid k \in \mathbf{Z}\}$ and set

$$(3.2) \qquad\qquad \varphi_s^\sigma(x) = \varphi_s(x)\mathbf{1}_{\{\sigma(x) \le |I_s|\}}.$$

The necessary result is

3.3. LEMMA. *Suppose that a collection of tiles $\mathbf{S}$ satisfies this combinatorial condition*:

$$(3.4) \qquad\qquad \omega_t \times I_t \cup \bigcup_{s \in \mathbf{T}_t} \omega_s \times I_s \qquad \text{are pairwise disjoint in } t.$$

*Suppose further that it satisfies the estimate*

$$(3.5) \qquad\qquad \sum_{s \in \mathbf{S}} |\langle f, \varphi_s \rangle|^2 \le K_0^2 \|f\|_2^2.$$

*Then for all $A, \mu \geq 1$, there are constants $C$ and $C_\mu$ for which the the collection $\mathbf{S}$ is a union of $\mathbf{S}^\sharp$ and $\mathbf{S}^\flat$, where the second collection $\mathbf{S}^\flat$ satisfies*

$$\left| \bigcup_{s \in \mathbf{S}^\flat} I_s \right| \leq C(A + \|N_{\mathbf{S}}\|_\infty)^{-\mu} \|N_{\mathbf{S}}\|_1.$$

*And for $\mathbf{S}^\sharp$, we have*

$$(3.6) \qquad \sum_{s \in \mathbf{S}^\sharp} |\langle f, \varphi_s^\sigma \rangle|^2 \leq C_\mu K_0^2 B^2 \|f\|_2^2,$$

$$(3.7) \qquad B = A^2 (\log A \|N_{\mathbf{S}}\|_\infty) \{ (\log A \|N_{\mathbf{S}}\|_\infty)^3 + A^{-\mu} \|N_{\mathbf{S}}\|_\infty^9 \}.$$

Section 4 of [7] gives specific combinatorial conditions under which (3.5) holds, with a quantitative estimate for $K_0$. Our estimate for $B$ is bigger than that for $K_0$ but the difference is harmless in the proof of the bounds for the bilinear maximal function. These combinatorial conditions are not strictly comparable to the condition we have imposed, (3.4). But both sets of conditions can be accommodated in the application of Lemma 3.3 to the proof of Theorem 2.14.

3.1. *Prerequisite maximal inequalities.* Maximal inequalities of two distinct types are needed to bound the operator in question. One is a simple variant of the Rademacher-Menschov theorem; the other is a variant of the Hardy-Littlewood maximal inequality due to Bourgain. Indeed this inequality is the critical ingredient of the pointwise ergodic theorem for arithmetic sets, [3]. To state it, consider basepoints $\lambda_\ell \in \mathbf{R}$, for $1 \leq \ell \leq L$, which are uniformly separated by $2^{-j_0}$. Thus, $|\lambda_\ell - \lambda_{\ell'}| \geq 2^{-j_0}$ for $\ell \neq \ell'$. Take $R_j$ to be a $2^{-j}$ neighborhood of these $L$ points, $R_j = \{\xi \mid \min_\ell |\xi - \lambda_\ell| \leq 2^{-j}\}$. We consider the Fourier restriction of $f$ to $R_j$,

$$(3.8) \qquad \Gamma_j f(x) = \int_{R_j} e(x\xi) \hat{f}(\xi) \, d\xi$$

and form a maximal operator from the $\Gamma_j$.

3.9. LEMMA. *We have the inequality*

$$\left\| \sup_{j \geq j_0} |\Gamma_j f| \right\|_2 \leq K(\log L)^3 \|f\|_2.$$

Notice that for $L = 1$ we recover the usual maximal function estimate, and that the significance of the result lies in the slow growth of the norm as a function of the number of base points. We refer the reader to Section 4 of [3] for the proof, which is a subtle integration of probabilistic and analytic techniques.

We need the following version of the Rademacher-Menschov theorem.

3.10. LEMMA. *Let $f_j$ be a sequence of functions on $L^2(X, \mu)$ with*

(3.11)
$$\Big\| \sum_{j=1}^{J} a_j f_j \Big\|_2 \le B \, \|a_j\|_\infty \, .$$

*Then we have the maximal inequality*

$$\Big\| \sup_{K < J} \Big| \sum_{j=1}^{K} f_j \Big| \Big\|_2 \le C(1 + \log J) B.$$

Such sequences $f_j$ are called "unconditionally convergent," and are well-known to enjoy many of the properties of orthogonal series. The most elegant result in this direction is one due to Ørno [8] which dilates the $f_j$ to an orthogonal series, but we cannot easily deduce our quantitative inequality from his. Our lemma is a close cousin of one due to G. Bennett, Theorem 2.5 of [1]; we do not include a proof because the lemma is easily obtained from the classical method of dyadic decomposition.

We also need the following corollary.

3.12. COROLLARY. *Let $f_n$ be a sequence of functions with*

$$\Big\| \sum_n a_n f_n \Big\|_2 \le B_0 \, \|a_n\|_\infty \, .$$

*Let $N_j$ be an increasing sequence of integers so that*

$$\Big\| \sup_{N_j \le N \le N_{j+1}} \Big| \sum_{n=N_j+1}^{N} f_n \Big| \Big\|_2 = B_j.$$

*Then we have the maximal inequality below, valid for all $J \ge 1$.*

$$\Big\| \sup_{N \le N_J} \Big| \sum_{n=1}^{N} f_n \Big| \Big\|_2 \le (1 + \log J) B_0 + \Big[ \sum_{j=1}^{J} B_j^2 \Big]^{1/2}.$$

*Proof.* Set $F_N = \sum_{n=1}^{N} f_n$, and let

$$F^j = \sup_{N_j \le N \le N_{j+1}} \Big| F_N - F_{N_j} \Big|.$$

Then

$$\Big\| \sup_{j \le J} \Big| F^j \Big| \Big\|_2^2 \le \sum_{j=1}^{J} \Big\| F^j \Big\|_2^2 \le \sum_{j=1}^{J} B_j^2.$$

And by Lemma 3.10,

$$\left\|\sup_{j \leq J}\left|F_{N_j}\right|\right\|_2 \leq (1 + \log J)B_0,$$

which completes the proof.                                                        □

3.2. *Combinatorics of grids.* The structure of grids enters into the proof as well. We will specifically need these lemmas, whose central concern is that of "enlarging" the space intervals in grids while maintaining the combinatorial structure of the grids.

3.13. LEMMA. *Let $A > 1$, and let $\mathbf{I}$ be a collection of intervals so that if $I \neq I' \in \mathbf{I}$ and $3/4\,|I| \leq |I'| \leq |I|$ then $I \cap I' = \emptyset$. Then there are intervals $I_A$ for $I \in \mathbf{I}$ so that $AI \subset I_A \subset (1+2^{-A})AI$ and $\{I_A \mid I \in \mathbf{I}\}$ is a union of $O(A^2)$ grids.*

*Proof.* We consider a special case. Let $\mathbf{I}$ be a collection of intervals which in addition satisfies for $I \neq I' \in \mathbf{I}$, the conditions that

$$\tfrac{3}{4}|I'| \leq |I| \leq |I'| \quad \text{implies} \quad \text{dist}(I, I') \geq 8A\,|I|\,, \text{ and that}$$
$$|I'| \leq \tfrac{3}{4}|I| \quad \text{implies} \quad |I'| \leq 2^{-A-3}\,|I|\,.$$

We then construct the $I_A$ so that $\{I_A \mid I \in \mathbf{I}\}$ is a grid. The special case proves the lemma, as the properties of $\mathbf{I}$ permit us to write it as a union of $O(A \log A) = O(A^2)$ subcollections which satisfy these last two properties.

Consider the graph with vertices $I \in \mathbf{I}$ and an edge from $I$ to $I'$ if $AI \cap AI' \neq \emptyset$. In this case we must have, e.g., $|I'| \leq 2^{-A-3}\,|I|$. For $I \in \mathbf{I}$ set $\mathbf{I}(I)$ to be those $I' \in \mathbf{I}$ which are connected to $I$ in the graph and $|I'| \leq |I|$. Then define $I_A = \bigcup_{I' \in \mathbf{I}(I)} AI'$. $I_A$ is an interval. Observe that $|I_A| \leq (1 + 2^{-A})\,|I|$. To see that $\{I_A \mid I \in \mathbf{I}\}$ is a grid, note that for $I \neq I' \in \mathbf{I}$ with $I_A \cap I'_A \neq \emptyset$ we must have e.g. $|I'| < 2^{-A}\,|I|$. But the intervals $I_A$ and $I'_A$ intersecting means that there must be an $I'' \in \mathbf{I}(I) \cap \mathbf{I}(I')$, and so each element of $\mathbf{I}(I')$ is connected to $I$, implying $I'_A \subset I_A$.                                □

We need a "separation lemma."

3.14. LEMMA. *Let $\mathbf{I}$ be a grid and let $N_{\mathbf{I}}(x) := \sum_{I \in \mathbf{I}} \mathbf{1}_I(x)$. Then for any $n \geq 2$ and $D > \|N_{\mathbf{I}}\|_\infty$, the grid can be split into two collections $\mathbf{I} = \mathbf{I}^\sharp \cup \mathbf{I}^\flat$ so that*

$$(3.15) \qquad \sum_{I \in \mathbf{I}^\sharp}\left(M\mathbf{1}_I\right)^2(x) \leq C_n D^3 \quad \text{for all } x,$$

$$(3.16) \qquad \left|\bigcup_{I \in \mathbf{I}^\flat} I\right| \leq 5D^{-n}\,\|N_{\mathbf{I}}\|_1\,.$$

*Proof.* We consider grids $\mathbf{I}$ for which

$$(3.17) \qquad I \subset_{\neq} I'\,;\; I, I' \in \mathbf{I} \quad \text{implies} \quad |I| \leq D^{-2n}\,|I'|\,.$$

For such a grid, we prove the lemma with $D^3$ in (3.15) replaced by $D^2$. As any grid is a union of $O(n \log D)$ subgrids satisfying (3.17) this is enough to prove the lemma.

The set which contains the "exceptional" $I$'s is determined in a two step procedure. In the first step, we set

$$I_{\partial} = \bigcup \Big\{ I' \in \mathbf{I} \mid I' \subset \{ x \mid \operatorname{dist}(x, \partial I) \leq D^{-n} |I| \} \Big\},$$

and take $E_1 = \bigcup_{I \in \mathbf{I}} I_{\partial}$. Obviously, $|E_1| \leq 4 D^{-n} \|N_{\mathbf{I}}\|_1$.

Then define a second exceptional set by letting

$$E_2 = \Big\{ x \mid \sum_{I \in \mathbf{I}} (M\mathbf{1}_I)^2(x) > C_n D^2 \Big\}.$$

Observe that the Fefferman-Stein maximal inequalities show that, for an appropriate $C_n$,

$$|E_2| \leq C_n^{-n} D^{-2n} \Big\| \sum_{I \in \mathbf{I}} (M\mathbf{1}_I)^2 \Big\|_n^n \leq D^{-2n} \Big\| \sum_{I \in \mathbf{I}} \mathbf{1}_I \Big\|_n^n \leq D^{-n} \|N_{\mathbf{I}}\|_1 \,.$$

We then take $\mathbf{I}^{\sharp} = \{ I \in \mathbf{I} \mid I \not\subset E_1 \cup E_2 \}$ and $\mathbf{I}^{\flat}$ is the complement of $\mathbf{I}^{\sharp}$. It is clear that (3.16) holds.

We verify that

$$\Big\| \sum_{I \in \mathbf{I}^{\sharp}} (M\mathbf{1}_I)^2 \Big\|_{\infty} \leq 17 K_n D^2.$$

And, as we argued at the beginning of the proof, from this we can conclude (3.15) holds. Indeed, assume that the inequality above does not hold; then there is an $x$ for which the sum above exceeds $17 K_n D^2$.

The assumption we prefer to argue from is that $x$ is an endpoint of one of the intervals $I$.

If $x$ is not an endpoint of an interval $I$, we remove from $\mathbf{I}^{\sharp}$ all those intervals $I \in \mathbf{I}^{\sharp}$ which strictly contain $x$. Then, the sum $\sum_{I \in \mathbf{I}^{\sharp}} (M\mathbf{1}_I)^2(x)$ still exceeds $16 K_n D^2$. The collection $\mathbf{I}^{\sharp}$ splits into those intervals $I \in \mathbf{I}^{\sharp}$ to the left of $x$ and those to the right. One of the two corresponding sums exceeds $8 K_n D^2$. If this is the case for those intervals to the right of $x$, we can increase the sum by moving $x$ to the right until it enters the union of the $I$'s. Hence for some endpoint $x$ of an $I_0 \in \mathbf{I}$ we have

$$(3.18) \qquad M(x) := \sum_{I \in \mathbf{I}^{\sharp}} (M\mathbf{1}_I)^2(x) \geq 8 C_n D^2.$$

This sum splits into

$$M_1(y) := \sum_{I \in \mathbf{I}^{\sharp}; |I| < |I_0|} (M\mathbf{1}_I)^2(y)$$

and $M_2(y) = M(y) - M_1(y)$. But, (3.17) and the removal of those $I \subset E_1$ imply that $M_1(x) \leq K$, where $K$ is an absolute constant. Thus, (3.18) implies that $M_2(x) \geq 7K_nD^2$. But this is a sum over $|I| \geq |I_0|$, and so we have

$$4 \inf_{y \in I_0} M\mathbf{1}_I(y)^2 \geq \sup_{y \in I_0} M\mathbf{1}_I(y)^2;$$

that is, we necessarily have $M_2(y) \geq C_n D^2$ for all $y \in I_0$. This contradicts the condition that $I_0 \not\subset E_2$ and so concludes the proof.  $\square$

3.3. *Main lemmas.* The argument to prove Lemma 3.3 is combinatorial and is achieved in several smaller steps. We first regularize the collection $\mathbf{S}$ by assuming that there is a grid $\{I_{s,A} \mid s \in \mathbf{S}\}$ so that $AI_s \subset I_{s,A} \subset (A + 2^{-A})I_s$ for $s \in \mathbf{S}$. We then prove the principal inequalities of this section without the factor of $A^2$ that appears on the right side of (3.7). This is enough to prove the estimates as stated, as Lemma 3.13 assures us.

To study the maximal inequality, we introduce a class of bounded operators. Consider the operators

$$T_\mathbf{S} f(x) = \sum_{s \in \mathbf{S}} \varepsilon_s \langle f, \varphi_{s,1} \rangle \varphi_{s,1}(x),$$

where $\varepsilon_s \in \{\pm 1\}$. Then, for all choices of sign

(3.19)                          $$\|T_\mathbf{S} f\|_2 \leq K_0 \mathcal{SQ}_\mathbf{S} f,$$

where $\mathcal{SQ}_\mathbf{S} f^2 = \sum_{s \in \mathbf{S}} |\langle f, \varphi_{s,1} \rangle|^2$. This is the dual form of the inequality (3.5). We introduce the notation of $\mathcal{SQ}_\mathbf{S} f$ because of our repeated applications of the Rademacher-Menschov theorem.

The class of maximal operators we study is

$$T_\mathbf{S}^{\max} f(x) = \sup_k \Big| \sum_{\substack{s \in \mathbf{S} \\ |I_s| \geq 2^k}} \varepsilon_s \langle f, \varphi_{s,1} \rangle \varphi_{s,1}(x) \Big|,$$

for once they are controlled, one can pass to a stopping time as in (3.2) and average over choices of signs to obtain the square function of Lemma 3.3.

A central theme is to decompose a collection $\mathbf{S}$ into a relatively large number of subcollections for which we can control the supremum. The "relatively large number" need not concern us because of the logarithmic term in the Rademacher-Menschov theorem, provided the subcollections fit together well. We codify these issues into

3.20. LEMMA. *Let $A, \mu \geq 1$. Suppose that $\mathbf{S}$ is a disjoint union of collections $\mathbf{S}_j$ for $1 \leq j \leq J$ so that for each $j$,*

(3.21)                          $$\Big\| T_{\mathbf{S}_j}^{\max} f \Big\|_2 \leq B \|f\|_2 + E \cdot \mathcal{SQ}_{\mathbf{S}_j} f.$$

*Moreover, the following combinatorial conditions relate the* $\mathbf{S}_j$. *Setting notation, each* $\mathbf{S}_j$ *is a union of* 1-*trees* $\mathbf{T}_t$ *with tops* $t \in \mathbf{S}_j^*$. *There are intervals* $\{\bar{I}_{j,v} \mid v \geq 1\}$ *such that*

(3.22) $\qquad\qquad \{\bar{I}_{j,v} \mid v \geq 1\}$ *are pairwise disjoint.*

(3.23) $\qquad\qquad$ *For each* $s \in \mathbf{S}_j$, *we have* $AI_s \subset \bar{I}_{j,v}$ *for some* $v$.

(3.24) $\qquad\qquad$ *If* $j < j'$ *then for all* $v'$, $\bar{I}_{j',v'} \subset \bar{I}_{j,v}$ *for some* $v$.

(3.25) $\qquad\qquad$ *For all* $j < j'$, $v \geq 1$, $s \in \mathbf{S}_j$ *with* $I_s \subset \bar{I}_{j,v}$

$\qquad\qquad\qquad$ *we have* $|I_s| \geq |I_{s'}|$ *if* $s' \in \mathbf{S}_{j'}$ *with* $I_{s'} \subset \bar{I}_{j,v}$.

*Then we have the following estimate, valid for all* $f \in L^2$ *of norm one.*

$$\left\| T_{\mathbf{S}}^{\max} f \right\|_2 \leq B \cdot \sqrt{J} + C_\mu A^{-\mu} \left\| N_{\mathbf{S}} \right\|_\infty \left\| \sum_{j,v} (M\mathbf{1}_{\bar{I}_{j,v}})^2 \right\|_\infty + (E + K_0 \log J) \mathcal{S}\mathcal{Q}_{\mathbf{S}} f.$$

In practice, we will have values like $\mu$ very large, $B = O(A^{-\mu} \|N_{\mathbf{S}}\|_\infty)$, $J = O(A \|N_{\mathbf{S}}\|_\infty^3)$ and $E = O(\log \|N_{\mathbf{S}}\|_\infty)$.

*Proof.* Fix $f \in L^2$ of norm one, and set for $s \in \mathbf{S}_j$, $I_s \subset \bar{I}_{j,v}$,

(3.26) $\qquad g_s(x) := \langle f, \varphi_{s,1} \rangle \varphi_{s,1}(x) \mathbf{1}_{\bar{I}_{j,v}}(x), \qquad h_s(x) := \langle f, \varphi_{s,1} \rangle \varphi_{s,1}(x) - g_s(x).$

The $h_s(x)$ are error terms. In particular, by (2.5), for $t \in \mathbf{S}_j^*$,

$$\sum_{\substack{s \in \mathbf{T}_t \\ I_s \subset \bar{I}_{j,v}}} |h_s(x)| \quad \leq \quad C \sum_{s \in \mathbf{T}_t} \sqrt{|I_s|} \, |\varphi_{s,1}(x)| \inf_{y \in I_s} Mf(y)$$

$$\leq \quad C_\mu A^{-\mu} \left( 1 + \frac{\operatorname{dist}(x, \bar{I}_{j,v})}{|\bar{I}_{j,v}|} \right)^{-3} M(Mf)(x),$$

so that

$$\sum_{s \in \mathbf{S}} |h_s(x)| \leq C_\mu A^{-\mu} M(Mf)(x) \left\| N_{\mathbf{S}} \right\|_\infty \left\| \sum_j (M\mathbf{1}_{\bar{I}_{j,v}})^2 \right\|_\infty.$$

It remains to bound the maximal sum over the $g_s(x)$, but here we observe that

(3.27) $\qquad \sup_k \left| \sum_{\substack{s \in \mathbf{S} \\ |I_s| \geq 2^k}} g_s(x) \right| \leq \sup_k \sup_K \left| \sum_{j=1}^{K-1} \sum_{s \in \mathbf{S}_j} g_s(x) + \sum_{\substack{s \in \mathbf{S}_K \\ |I_s| \geq 2^k}} g_s(x) \right|.$

Indeed, this is an easy consequence of (3.25). Fix $x$ and $k$. Let $\bar{s}$ be a tile in $\mathbf{S}$ so that $x \in \operatorname{supp}(g_s)$ and $|I_{\bar{s}}| \geq 2^k$, but $|I_{\bar{s}}|$ is in addition minimal. Then $\bar{s} \in \mathbf{S}_K$ for some $K$, and we observe that

$$\sum_{\substack{s \in \mathbf{S} \\ |I_s| \geq 2^k}} g_s(x) = \sum_{j=1}^{K-1} \sum_{s \in \mathbf{S}_j} g_s(x) + \sum_{\substack{s \in \mathbf{S}_K \\ |I_s| \geq 2^k}} g_s(x).$$

By construction, there can be no $s' \in \bigcup_{j>K} \mathbf{S}_{j'}$ with $g_{s'}(x) \neq 0$ and $|I_{s'}| > 2^k$. So, suppose that $s \in \mathbf{S}_j$ with $j < K$ and $g_s(x) \neq 0$. Now, $I_{\bar{s}}, I_s \subset \bar{I}_{j,v}$ for some $v$. Then, (3.25) implies that $|I_s| \geq |I_{\bar{s}}| \geq 2^k$. This proves the equality above.

Finally, (3.21) and Corollary 3.12 imply that the $L^2$ norm of the right-hand side of (3.27) is bounded by

$$B\sqrt{J} + C_\mu A^{-\mu} \|N_\mathbf{S}\|_\infty \left\|\sum_{j,v} (M\mathbf{1}_{\bar{I}_{j,v}})^2\right\|_\infty + (E + K_0 \log J)\mathcal{S}\mathcal{Q}_\mathbf{S} f.$$

This completes the proof.                                                      □

We now concern ourselves with trees that are "trivial" in that they consist of only a bounded number of layers.

**3.28. Lemma.** *Let* $\mathbf{S}$ *be as in Lemma* 3.3 *with the additional assumption that for all for all* $t \in \mathbf{S}^*$

$$(3.29) \qquad\qquad \sum_{s \in \mathbf{T}_t} \mathbf{1}_{I_s}(x) \leq A_0 \quad \textit{for all } x.$$

*Then if* $f$ *has* $L^2$ *norm one,*

$$
\begin{aligned}
(3.30) \qquad \|T_\mathbf{S}^{\max} f\|_2 \quad \leq \quad & C_\mu A^{-\mu} A_0^2 \left\|\sum_{s \in \mathbf{S}} (M\mathbf{1}_{I_s})^2\right\|_\infty^2 \\
& + CK_0\left(1 + \log A_0 \left\|\sum_{t \in \mathbf{S}^*} \mathbf{1}_{2AI_t}\right\|_\infty\right)\mathcal{S}\mathcal{Q}_\mathbf{S} f.
\end{aligned}
$$

*Proof.* The following decomposition of $\mathbf{S}$ is central to the proof of the lemma. Recall that we have a grid $\{I_{s,A} \mid s \in \mathbf{S}\}$ with $AI_s \subset I_{s,A} \subset 2AI_s$ for all $s$. Then let $\mathbf{S}_1$ be those $s \in \mathbf{S}$ for which $I_{s,A}$ is maximal. Remove these tiles from $\mathbf{S}$ and repeat to define $\mathbf{S}_2$ and so on. This procedure must stop in $J = O(A_0 \|\sum_{t \in \mathbf{S}^*} \mathbf{1}_{2AI_t}\|_\infty)$ steps.

It is our intention to apply Lemma 3.20 to this decomposition of $\mathbf{S}$. To this end, we take the collections of intervals $\{\bar{I}_{j,v} \mid v \geq 1\}$ in that lemma to be any enumeration of the intervals $\{I_{s,A} \mid s \in \mathbf{S}_j\}$. The four conditions (3.22)–(3.25) follow immediately. So by Lemma 3.20 it suffices to bound a single $T_{\mathbf{S}_j}^{\max}$.

But this is a relatively easy matter, for the rectangles in $\mathbf{S}_j$ are maximal by construction. We define the functions $g_s$ as in (3.26), namely,

$$
\begin{aligned}
g_s(x) \quad &:= \quad \langle f, \varphi_{s,1}\rangle \varphi_{s,1}(x)\mathbf{1}_{I_{s,A}}(x), \\
h_s(x) \quad &:= \quad \langle f, \varphi_{s,1}\rangle \varphi_{s,1}(x) - g_s(x).
\end{aligned}
$$

Hence, for any $x \in \mathbf{R}$ and $s, s' \in \mathbf{S}_j$ if $g_s(x)g_{s'}(x) \neq 0$ then $I_s = I_{s'}$. Thus, if $f$ has norm one,

$$
\begin{aligned}
\left\|T^{\max}_{\mathbf{S}_j} f\right\|_2 &\leq \left\|\sum_{s \in \mathbf{S}_j} |h_s|\right\|_2 + \left\|\sum_{s \in \mathbf{S}_j} g_s\right\|_2 \\
&\leq C_\mu A^{-\mu} \left\|\sum_{s \in \mathbf{S}_j} (M\mathbf{1}_{I_s})^2\right\|_\infty + \left\|T_{\mathbf{S}_j} f\right\|_2 \\
&\leq C_\mu A^{-\mu} \left\|\sum_{s \in \mathbf{S}} (M\mathbf{1}_{I_s})^2\right\|_\infty + K_0 \mathcal{S}\mathcal{Q}_{\mathbf{S}} f.
\end{aligned}
$$

Therefore, the lemma follows. $\qquad\blacksquare$

We turn to the case of nontrivial trees. Indeed, this is the crucial issue. Here some technical issues become important, issues that arise from our formulation of the Rademacher-Menschov theorem and the form of Bourgain's lemma.

We consider $T^{\max}_{\mathbf{T}} f$, where $\mathbf{T}$ is a 1-tree and identify it as the maximum of Fourier projections applied to $T_{\mathbf{T}} f$. We separate scales in a manner dictated by the parameters $a$ and $K$ that appear in (2.4) and (2.13). We shall assume that for all tiles $s$

$$
(3.31) \qquad 2^{\alpha k + \beta} \leq |I_s| \leq \tfrac{4}{3} 2^{\alpha k + \beta} \qquad \text{for some } k.
$$

Here, $\alpha = \alpha(K, a)$ is a sufficiently large integer and $\beta \in \{0, 1, \dots, \alpha - 1\}$ is fixed.

Consider a 1-tree $\mathbf{T}_t$ with top $t$. Let $\lambda_t = a(c(\omega_t))$, where $c(\omega_t)$ is the center of $\omega_t$ and $a$ is the affine function as in (2.4). Then the Fourier transform of $\varphi_{s,1}$ is supported on $v_{s,1} := \frac{3}{4}\{a(\xi) \mid \xi \in \omega_{s,1}\}$, and $\omega_t \cap \omega_{s,1}$ is empty. Let $a(x) = m_a x + b$. Hence from (2.3) and (2.13)

$$
\begin{aligned}
\tilde{K} |m_a| |I_s|^{-1} &\leq \tfrac{1}{4} |m_a \omega_{s,1}| \\
&\leq \inf\{|\xi - \lambda_t| \mid \xi \in v_{s,1}\} \\
&\leq \sup\{|\xi - \lambda_t| \mid \xi \in v_{s,1}\} \\
&\leq |m_a| |\omega_s| \\
&\leq K |m_a| |I_s|^{-1}.
\end{aligned}
$$

The constants $K$ and $\tilde{K}$ depend only on (2.3). Namely,

$$
\tilde{K} := \tfrac{1}{4} \inf_s |s| < K := \sup_s |s|.
$$

While these constants can depend on the choice of model sum, (2.3) shows that the ratio is independent of the model sum. We take $\alpha$ to be an integer with $2^\alpha > K/\tilde{K}$ and in addition choose an integer $\gamma$ with $2^\gamma > 2K\,|m_a|$. With these choices, for integers $k$

$$(3.32) \qquad \sum_{\substack{s \in \mathbf{T}_t \\ |I_s| \geq 2^k}} \varepsilon_s \langle f, \varphi_{s1} \rangle \varphi_{s1}(x) = \int_{J(k)} e(x\xi) \widehat{T_{\mathbf{T}_t} f}(\xi)\, d\xi,$$

where $J(k) = (\lambda_t - 2^{-k+\gamma}, \lambda_t + 2^{-k+\gamma})$. This is designed to fit well with the Fourier restriction operators of (3.8).

With this observation behind us, we state a lemma; this admits an especially direct appeal to Bourgain's lemma; this requires a nice structure on $\mathbf{S}$, which we formalize with this definition.

**3.33.** *Definition.* Call $\mathbf{S}$ *uniform* if it satisfies the hypotheses of Lemma 3.3, $\bigcap_{t \in \mathbf{S}^*} I_t \neq \emptyset$ and

$$\text{dist}(\omega_t, \omega_{t'}) \geq \beta \max_{s \in \mathbf{S}} |\omega_s|, \quad t \neq t' \in \mathbf{S}^*.$$

Here the quantity $\beta := \inf |\omega_s 1|\,/\,|\omega_s|$, with the infimum taken over all tiles $s$. By (2.3) and (2.13) this is a positive constant.

**3.34.** LEMMA. *Let $\mu, A \geq 1$. If $\mathbf{S}$ is uniform and $D \geq A + \|N_{\mathbf{S}}\|_\infty$, then $\mathbf{S}$ is the union of $\mathbf{S}^\sharp$ and $\mathbf{S}^\flat$, where*

$$(3.35) \qquad \left| \bigcup_{s \in \mathbf{S}^\flat} I_s \right| \leq C_\mu D^{-\mu} \|N_{\mathbf{S}}\|_1.$$

*And for $\mathbf{S}^\sharp$ and all $f \in L^2$ of norm one,*

$$\left\| T_{\mathbf{S}^\sharp}^{\max} f \right\|_2 \leq C_\mu A^{-\mu} D^6 + K_0 \psi_3(AD) \mathcal{S} \mathcal{Q}_{\mathbf{S}} f.$$

*In this last display, we use the notation*

$$(3.36) \qquad \psi_j(d) := C_\mu (1 + \log d)^j, \quad d,\ j \geq 1.$$

*Proof.* The essence of the matter is (3.32) and Lemma 3.9. In order to use these facts, we need to remove the bottom from $\mathbf{S}$. To this end, take $\mathbf{S}_{\text{bott}}$ to consist of all $s \in \mathbf{S}$ for which

$$|I_s| < (\beta\,|m_a|\,\tilde{K})^{-1} 2^\gamma \inf_{s' \in \mathbf{S}} |I_{s'}|.$$

It follows that Lemma 3.28 applies to $\mathbf{S}_{\text{bott}}$, with $A_0 = O(\gamma + |\log \alpha| + |\log \beta|)$. Each of these last three quantities can be regarded as fixed; that is, $A_0$ can be taken to be a constant. But in addition we can apply Lemma 3.14 to the grid $\{I_s \mid s \in \mathbf{S}_{\text{bott}}\}$ to select a subset $\mathbf{I}^\flat$ with measure as in (3.35), so that for the collection $\mathbf{S}_0^\sharp := \{s \in \mathbf{S}_{\text{bott}} \mid I_s \notin \mathbf{I}^\flat\}$ we have

$$\left\| \sum_{s \in \mathbf{S}_0^\sharp} (M\mathbf{1}_{I_s})^2 \right\|_\infty \leq C_\mu D^3.$$

Hence, for all $f$ of $L^2$ norm one, we have the estimate

$$\left\| T^{\max}_{\mathbf{S}^{\sharp}_0} f \right\|_2 \le C_\mu A^{-\mu} D^6 + C K_0 \psi_1(AD) \mathcal{S} \mathcal{Q}_{\mathbf{S}^{\sharp}_0} f.$$

This follows from Lemma 3.28.

As for $\mathbf{S}^{\sharp} = \mathbf{S} \backslash \mathbf{S}_{\text{bott}}$, there is no need delete additional tiles. Now we may have deleted all of a tree $\mathbf{T}_t$ in $\mathbf{S}_{\text{bott}}$. Thus, delete from $\mathbf{S}^*$ any top with $\mathbf{T}_t \subset \mathbf{S}_{\text{bott}}$. For $t \ne t' \in \mathbf{S}^*$, recall that $\lambda_t = a(c(\omega_t))$, where $a$ is the affine function in (2.4). Then

$$
\begin{aligned}
|\lambda_t - \lambda_{t'}| &\ge |m_a| \operatorname{dist}(\omega_t, \omega_{t'}) \\
&\ge \beta |m_a| \sup_{s \in \mathbf{S}} |\omega_s| \\
&\ge \tilde{K} |m_a| \beta \sup_{s \in \mathbf{S}} |I_s|^{-1} \\
&\ge 2^\gamma \sup_{s \in \mathbf{S}^{\sharp}} |I_s|^{-1} .
\end{aligned}
$$

Hence, we have the equality below, following from (3.32).

$$\sum_{\substack{s \in \mathbf{S}^{\sharp} \\ |I_s| \ge 2^k}} \varepsilon_s \langle f, \varphi_{s,1} \rangle \varphi_{s,1}(x) = \int_{R(k)} e(x\xi) \widehat{T_{\mathbf{S}^{\sharp}} f}(\xi) \, d\xi,$$

where $R(k) = \bigcup_{t \in \mathbf{S}^*} (\lambda_t - 2^{-k+\gamma}, \lambda_t + 2^{-k+\gamma})$. Therefore, Lemma 3.9 implies that

$$\left\| T^{\max}_{\mathbf{S}^{\sharp}} f \right\|_2 \le C(\log \|N_{\mathbf{S}}\|_\infty)^3 \|T_{\mathbf{S}^{\sharp}} f\|_2,$$

and then (3.5) finishes the proof. $\qquad\square$

We enlarge the class of tiles to which the previous estimate applies in three distinct steps.

**3.37.** *Definition.* Call $\mathbf{S}$ *separated* if it satisfies the hypotheses of Lemma 3.3 and there is a sequence of intervals $\bar{I}_v$ with $\{A\bar{I}_v \mid v \ge 1\}$ pairwise disjoint, and $\mathbf{S}$ is a union of collections $\mathbf{S}_v := \{s \in \mathbf{S} \mid I_s \subset \bar{I}_v\}$, each of which is uniform.

**3.38.** **Lemma.** *Let $A, \mu \ge 1$. If $\mathbf{S}$ is separated and $D \ge A + \|N_{\mathbf{S}}\|_\infty$, then $\mathbf{S} = \mathbf{S}^{\sharp} \cup \mathbf{S}^{\flat}$, with*

$$(3.39) \qquad \left| \bigcup_{s \in \mathbf{S}^{\flat}} I_s \right| \le C_\mu D^{-\mu} \left\{ \|N_{\mathbf{S}}\|_1 + \sum_v |\bar{I}_v| \right\},$$

*and for the second collection of tiles $\mathbf{S}^{\sharp}$ and for all $f$ of $L^2$ norm one,*

$$\left\| T^{\max}_{\mathbf{S}^{\sharp}} f \right\|_2 \le C_\mu A^{-\mu} D^6 + K_0 \psi_3(AD) \mathcal{S} \mathcal{Q}_{\mathbf{S}} f.$$

*Proof.* There is an initial contribution to the collection $\mathbf{S}^{\flat}$ to make. The collection $\{A\bar{I}_v\}$ is trivially a grid, since they are disjoint. So by Lemma 3.14, the collection of intervals $\{\bar{I}_v\}$ is the union of $\{\bar{I}_v \mid v \in V^\sharp\}$ and $\{\bar{I}_v \mid v \in V^\flat\}$ with

$$(3.40) \qquad \left\| \sum_{v \in V^\sharp} (M\mathbf{1}_{\bar{I}_v})^2 \right\|_\infty \le D^3,$$

$$\left| \bigcup_{v \in V^\flat} \bar{I}_v \right| \le C_\mu D^{-\mu} \sum_v |\bar{I}_v|.$$

We take our initial contribution to $\mathbf{S}^\flat$ to be $\mathbf{S}_0^\flat := \{s \in \mathbf{S} \mid I_s \subset \bar{I}_v, \ v \in V^\flat\}$. We do not further consider these tiles. Thus, we remove them from $\mathbf{S}$.

The inequality is obvious for each $\mathbf{S}_v := \{s \in \mathbf{S} \mid I_s \subset \bar{I}_v\}$. We apply Lemma 3.34 to $\mathbf{S}_v$. Hence $\mathbf{S}_v = \mathbf{S}_v^\sharp \cup \mathbf{S}_v^\flat$, with

$$\left| \bigcup_{s \in \mathbf{S}_v^\flat} I_s \right| \le C_\mu D^{-\mu} \|N_{\mathbf{S}_v}\|_1.$$

And for $\mathbf{S}_v^\sharp$, we have

$$\left\| T_{\mathbf{S}_v^\sharp}^{\max} f \right\|_2 \le C_\mu A^{-\mu} D^6 + \psi_3(AD) \mathcal{S} \mathcal{Q}_{\mathbf{S}} f,$$

for all $f$ of $L^2$ norm one.

But it is clear that $T_{\mathbf{S}_v^\sharp}$ is nearly supported on $A\bar{I}_v$. Namely, for all $x \in \mathbf{R}$,

$$\mathbf{1}_{(A\bar{I}_v)^c}(x) \sum_{s \in \mathbf{S}_v^\sharp} |\langle f, \varphi_{s,1} \rangle \varphi_{s,1}(x)| \le C_\mu A^{-\mu} \|N_{\mathbf{S}}\|_\infty (M\mathbf{1}_{\bar{I}_v})^2(x) Mf(x).$$

Since (3.40) holds, this proves the lemma. $\qquad \square$

3.41. *Definition.* Call $\mathbf{S}$ *normal* if it satisfies the hypotheses of Lemma 3.3 and there are intervals $\{\bar{I}_v \mid v \ge 1\}$ for which $\{A\bar{I}_v \mid v \ge 1\}$ are pairwise disjoint, and for all $t \in \mathbf{S}^*$ and all $s \in \mathbf{T}_t$, there is a $v \ge 1$ so that $I_t \supset \bar{I}_v \supset I_s$.

3.42. Lemma. *If $\mathbf{S}$ is normal and $D \ge A + \|N_{\mathbf{S}}\|_\infty$, then $\mathbf{S} = \mathbf{S}^\sharp \cup \mathbf{S}^\flat$, with*

$$(3.43) \qquad \left| \bigcup_{s \in \mathbf{S}^\flat} I_s \right| \le C_\mu D^{-\mu} \|N_{\mathbf{S}}\|_1,$$

*and for the second collection of tiles $\mathbf{S}^\sharp$ and for all $f$ of $L^2$ norm one,*

$$(3.44) \qquad \left\| T_{\mathbf{S}}^{\max} f \right\|_2 \le C_\mu A^{-\mu} D^8 + K_0 \psi_3(AD) \mathcal{S} \mathcal{Q}_{\mathbf{S}} f.$$

*Proof.* Just as in the previous proof, we make an initial contribution to the collection $\mathbf{S}^\flat$. The collection of intervals $\{\bar{I}_v\}$ is the union of $\{\bar{I}_v \mid v \in V^\sharp\}$ and $\{\bar{I}_v \mid v \in V^\flat\}$ with

$$(3.45) \qquad \left\| \sum_{v \in V^\sharp} (M\mathbf{1}_{\bar{I}_v})^2 \right\|_\infty \leq D^3,$$

$$\left| \bigcup_{v \in V^\flat} \bar{I}_v \right| \leq C_\mu D^{-\mu} \sum_v |\bar{I}_v| \leq C_\mu D^{-\mu} \|N_\mathbf{S}\|_1 .$$

Define $\mathbf{S}_0^\flat := \{s \in \mathbf{S} \mid I_s \subset \bar{I}_v, \ v \in V^\flat\}$. We remove them from $\mathbf{S}$.

We apply Lemma 3.20, where each $\mathbf{S}_j$ is separated and the number of $\mathbf{S}_j$ is $J \leq \|N_\mathbf{S}\|_\infty$.

For each $v$, consider the tops $t(v,j) \in \mathbf{S}^*$ with $1 \leq j \leq J_v \leq \|N_\mathbf{S}\|_\infty$ and $\bar{I}_v \subset I_{t(v,j)}$. Let

$$\delta(v,j) := \inf\{|I_s| \mid s \in \mathbf{T}_{t(v,j)}, \ I_s \subset \bar{I}_v\}.$$

We may assume that the $\delta(v,j)$ decrease in $j$. Then take

$$\mathbf{S}_{v,j} := \Big\{ s \in \mathbf{T}_{t(v,j)} \mid \delta(v,j) \geq |I_s| > \delta(v,j+1) \Big\},$$

and $\mathbf{S}_j = \bigcup_v \mathbf{S}_{v,j}$. To apply Lemma 3.20, we take the collections $\{\bar{I}_{j,v} \mid v \geq 1\}$ to be $\{\bar{I}_v\}$ as specified in the definition of normal. By construction, the $\mathbf{S}_j$ satisfy (3.22)–(3.25).

Moreover, each $\mathbf{S}_{j,v}$ is uniform. For if $t \neq t' \in \mathbf{S}^*$ and there are $s \in \mathbf{T}_t \cap \mathbf{S}_{j,v}$ and $s' \in \mathbf{T}_{t'} \cap \mathbf{S}_{j,v}$ then by construction there are $\bar{s} \in \mathbf{T}_t$ with $|I_{\bar{s}}| \leq \delta(j+1,v)$ and likewise for $\bar{s}'$. But the crucial combinatorial property (3.4) then enters in. It implies that $\omega_{\bar{s}}$ and $\omega_{\bar{s}'}$ are disjoint. By (2.11), either $\omega_{\bar{s}}1$ or $\omega_{\bar{s}'}1$ lies between $\omega_t$ and $\omega_{t'}$. Hence

$$\mathrm{dist}(\omega_t, \omega_{t'}) \geq \beta \inf\{|\omega_{\bar{s}}|, |\omega_{\bar{s}'}|\} \geq \beta \sup_{s \in \mathbf{S}_{j,v}} |\omega_s|,$$

proving uniformity. ($\beta$ is defined in Definition 3.33.)

Hence each $\mathbf{S}_j$ is separated. Apply Lemma 3.38 with $D = A + \|N_\mathbf{S}\|_\infty$. This splits $\mathbf{S}_j$ into $\mathbf{S}_j^\sharp$ and $\mathbf{S}_j^\flat$. In the estimate for $\mathbf{S}_j^\flat$, (3.39), we have $\sum_v |\bar{I}_v| \leq C_\mu D^{-\mu} \|N_\mathbf{S}\|_1$ and we use (3.39) at most $J \leq \|N_\mathbf{S}\|_\infty \leq D$ times.

In addition, $\sum_{v \in V^\sharp} (M\mathbf{1}_{\bar{I}_v})^2$ is uniformly controlled, (3.45), by the first step in the argument. The lemma then follows from Lemma 3.20. $\qquad \blacksquare$

*Proof of Lemma* 3.3. The idea of the proof is to write $\mathbf{S}$ as a union of at most $\|N_\mathbf{S}\|_\infty$ normal families $\mathbf{S}_j$. There is extra difficulty in that we cannot arrange an application of Lemma 3.20 for these collections of tiles, so an additional argument must be made to conclude the lemma.

Take $D = A + \|N_\mathbf{S}\|_\infty$, $f$ of $L^2$ norm one and assume that $\mu \geq 400$.

We can assume for all $s, s' \in \mathbf{S}$ that $|I_s| < \frac{3}{4}|I_{s'}|$ implies $D^{4\mu}|I_s| < |I_{s'}|$. (To explain the role of $\frac{3}{4}$, see the definition of a grid, Definition 2.7.) And prove the lemma without a power of $\log(A \|N_\mathbf{S}\|_\infty)$ in the expression for $B$ in (3.7). This is clearly sufficient to conclude the lemma as stated.

Several distinct contributions to the set $\mathbf{S}^\flat$ must be made. The first contribution, defined now, plays a role in the final argument of this proof. The set

$$F := \bigcup_{t \in \mathbf{S}^*} \{ x \mid \operatorname{dist}(x, \partial I_t) < D^{-\mu} \, |I_t| \}$$

has measure $|F| \leq 2D^{-\mu} \, \|N_{\mathbf{S}}\|_1$. Take $\mathbf{S}^\flat_{-1}$ to be $\{ s \in \mathbf{S} \mid I_s \subset F \}$, and remove these tiles from $\mathbf{S}$.

We make another contribution to $\mathbf{S}^\flat$. By Lemma 3.14 applied to $\{ I_t \mid t \in \mathbf{S}^* \}$, the collection of tops $\mathbf{S}^*$ is the union of $\mathbf{S}^{*\sharp}$ and $\mathbf{S}^{*\flat}$ with

$$(3.46) \qquad \left| \bigcup_{t \in \mathbf{S}^{*\flat}} I_t \right| \;\leq\; C_\mu D^{-\mu} \, \|N_{\mathbf{S}}\|_1 \, ,$$

$$\left\| \sum_{s \in \mathbf{S}^{*\sharp}} (M \mathbf{1}_{I_t})^2 \right\|_\infty \;\leq\; C_\mu \, \|N_{\mathbf{S}}\|_\infty^3 \, .$$

We then take $\mathbf{S}^\flat_0 := \{ s \in \mathbf{S} \mid I_s \subset I_t, \ t \in \mathbf{S}^{*\flat} \}$ to be our second contribution to $\mathbf{S}^\flat$. These tiles are removed from $\mathbf{S}$.

Recall that we began this section with the assumption that there is a grid $\mathcal{G} = \{ I_{t,A} \mid t \in \mathbf{S}^* \}$ with $AI_t \subset I_{t,A}$ for all $t \in \mathbf{S}^*$. Let $\mathcal{G}_1$ consist of the maximal elements of $\mathcal{G}$, and inductively define $\mathcal{G}_{j+1}$ to be the maximal elements of $\mathcal{G} \backslash \bigcup_{\ell=1}^j \mathcal{G}_\ell$. Clearly this procedure stops in at most $J \leq \|N_{\mathbf{S}}\|_\infty$ steps.

Let $\mathbf{S}_J := \{ s \in \mathbf{S} \mid I_s \subset I_t, \ I_{t,A} \in \mathcal{G}_J \text{ for some t} \}$, and by reverse induction define

$$\mathbf{S}_{j-1} = \Big\{ s \in \mathbf{S} \backslash \bigcup_{\ell=j}^J \mathbf{S}_\ell \mid I_s \subset I_t, \ I_{t,A} \in \mathcal{G}_{j-1} \text{ for some t} \Big\}.$$

Each $\mathbf{S}_j$ is normal, as is easy to see. Thus Lemma 3.42 applies to each $\mathbf{S}_j$, yielding $\mathbf{S}_j = \mathbf{S}_j^\sharp \cup \mathbf{S}_j^\flat$ with $\mathbf{S}_j^\flat$ as in (3.43) and $\mathbf{S}_j^\sharp$ as in (3.44). Recall that $1 \leq j \leq J \leq D$.

Then take $\mathbf{S}^\flat = \bigcup_{j=-1}^J \mathbf{S}_j^\flat$. It is clear that the desired estimate holds on $\mathbf{S}^\flat$, and we remove all of these tiles from $\mathbf{S}$ and each $\mathbf{S}_j$.

The remaining difficulty is that Lemma 3.20 does not directly apply to the $\mathbf{S}_j$. And indeed, we must adopt a more sophisticated truncation argument to complete the proof.

We remove the top from $\mathbf{S}$. Let $\mathbf{S}_j^{\text{tops}}$ be those $t \in \mathbf{S}_j$ for which $I_{t,A} \in \mathcal{G}_j$. Set $\mathbf{S}^{\text{tops}} = \bigcup_j \mathbf{S}_j^{\text{tops}}$. By Lemma 3.28, applied with $A_0 = 1$, and (3.46), we have

$$\left\| T_{\mathbf{S}^{\text{tops}}}^{\max} f \right\|_2 \leq C_\mu A^{-\mu} D^6 + \psi_1(AD) \mathcal{S} \mathcal{Q}_{\mathbf{S}^{\text{tops}}} f.$$

But in addition, the collections $\mathbf{S}_j^{\text{tops}}$ satisfy the hypotheses of Lemma 3.20. Therefore, these tiles satisfy the conclusion of the lemma, and we can remove them from $\mathbf{S}$.

Let $\mathcal{I}_j$ consist of the intervals $I_t$ with $I_{t,A} \in \mathcal{G}_j$. Recall from the first step of the proof, and the removal of the tops, that for each $I \in \mathcal{I}_j$, if $s \in \mathbf{S}_j$ with $|I_s| < \frac{3}{4}|I|$, then $D^{4\mu}|I_s| < |I|$, and

$$(3.47) \qquad \mathrm{dist}(I_s, \partial I) > D^{-\mu}|I|.$$

As well, $\mu > 400$.

For $s \in \mathbf{S} = \bigcup \mathbf{S}_j$, set $f_s = \langle f, \varphi_{s,1} \rangle \varphi_{s,1}$. We will write $f_s = a_s + b_s + c_s$ where $a_s$ is the dominant term. To define these terms, suppose that $s \in \mathbf{S}_j$, $I_s \subset I$ with $I \in \mathcal{I}_j$. Then define

$$b_s(x) := f_s(x)\mathbf{1}_{I^c}(x).$$

To define $c_s$, let $\mathcal{I}(j,s)$ be those $I' \in \mathcal{I}_{j'}$ with $j' > j$, $I' \subset I$ and $|I_s| < \frac{3}{4}|I'|$. Then define

$$(3.48) \qquad c_s(x) := f_s(x)\mathbf{1}_{E_s}(x), \qquad E_s := \bigcup \{I' \mid I' \in \mathcal{I}(j,s)\}.$$

Finally, $a_s = f_s - b_s - c_s$.

The principal point of this construction is that

$$(3.49) \qquad \sup_n \Big| \sum_{\substack{s \in \mathbf{S} \\ |I_s| \geq 2^n}} a_s(x) \Big| \leq \sup_k \sup_n \Big| \sum_{j=1}^{k-1} \sum_{s \in \mathbf{S}_j} a_s(x) + \sum_{\substack{s \in \mathbf{S}_k \\ |I_s| \geq 2^n}} a_s(x) \Big|.$$

We will return to this momentarily.

In addition, the terms $b_s$ and $c_s$ are error terms.

$$(3.50) \qquad \Big\| \sum_{s \in \mathbf{S}} |b_s| + |c_s| \Big\|_2 \leq C_\mu D^{-\mu}.$$

We prove this below.

But assuming these last two claims, we see that by construction and normality for the $\mathbf{S}_j$, we have for each $1 \leq j \leq D = \|N_{\mathbf{S}}\|_\infty$,

$$\Big\| \sup_n \Big| \sum_{\substack{s \in \mathbf{S}_j \\ |I_s| \geq 2^n}} a_s(x) \Big| \Big\|_2 \leq C_\mu A^{-\mu} D^8 + K_0 \psi_3(AD) \mathcal{S} \mathcal{Q}_{\mathbf{S}_j} f.$$

Therefore, by (3.49) and Corollary 3.12, we have

$$\Big\| \sup_n \Big| \sum_{\substack{s \in \mathbf{S} \\ |I_s| \geq 2^n}} a_s(x) \Big| \Big\|_2 \leq C_\mu A^{-\mu} D^9 + K_0 \psi_3(AD) \mathcal{S} \mathcal{Q}_{\mathbf{S}} f,$$

which is the estimate we need to conclude the proof.

We now turn to the principal inequality, (3.49). Fix an integer $n$ and $x \in \mathbf{R}$. Let $j'$ be the maximal integer such that $a_{s'}(x) \neq 0$ for some $s' \in \mathbf{S}_{j'}$ with $|I_{s'}| \geq 2^n$. Then our task is to show that for any $j < j'$ and $s \in \mathbf{S}_j$ with $|I_s| < 2^n$, we have $a_s(x) = 0$. But, choosing $I' \in \mathcal{I}_{j'}$ with $I_{s'} \subset I'$, and $I \in \mathcal{I}_j$ with $I_s \subset I$, we are free to assume $x \in I' \subset I$, else there is nothing to prove. We would then have

$$|I_s| < 2^n \leq |I_{s'}| < \tfrac{3}{4} |I'| .$$

Then it follows from the definition (3.48) that $I' \subset E_s$. But $x \in I_{s'} \subset I'$ by assumption, so $a_s(x) = 0$. (3.49) is proved.

The remaining element of the proof is the verification of (3.50). Let us treat the terms $b_s$ first. By (3.47), we have $D^\mu I_s \subset I$ for $I \in \mathcal{I}_j$, $s \in \mathbf{S}_j$ and $I_s \subset I$. Hence, for all $x \in \mathcal{R}$,

$$
\begin{aligned}
\mathbf{1}_{I^c}(x) \sum_{\substack{s \in \mathbf{S}_j \\ I_s \subset I}} |b_s(x)| &\leq C_\mu D^{-2\mu} \sum_{\substack{s \in \mathbf{S}_j \\ I_s \subset I}} \left(1 + \frac{\operatorname{dist}(x,I)}{|I_s|}\right)^{-20} \inf_{y \in I_s} Mf(y) \\
&\leq C_\mu D^{-2\mu} \|N_{\mathbf{S}}\|_\infty \sum_{I \in \mathcal{I}_j} (M\mathbf{1}_I)^2 MMf(x).
\end{aligned}
$$

And, by (3.46), the $L^2$ norm of this last expression is at most $C_\mu D^{-2\mu+3}$. But $\mu \geq 400$, so (3.50) is verified for the sum over $b_s$.

To estimate the sum over $c_s$, recall the definition (3.48) to see that

$$
\begin{aligned}
\sum_{s \in \mathbf{S}} |c_s(x)| &\leq \sum_{s \in \mathbf{S}} \sum_{I' \in \mathcal{I}(j,s)} |f_s(x)\mathbf{1}_{I'}(x)| \\
&\leq \sum_{j'=1}^{J} \sum_{I' \in \mathcal{I}_{j'}} \mathbf{1}_{I'}(x) \sum_{s \in \mathcal{J}(j',I')} |f_s(x)| .
\end{aligned}
$$

In this last sum, the collection of tiles $\mathcal{J}(j',I')$ consists of those tiles $s \in \bigcup_{j=1}^{j'-1} \mathbf{S}_j$ for which $I' \subset E_s$. Now, for such $s$, we necessarily have $D^{4\mu} |I_s| < |I'|$ and $\operatorname{dist}(I_s, \partial I') \geq D^{-\mu} |I'|$. Fixing $I'$ above, we have

$$
\begin{aligned}
\mathbf{1}_{I'}(x) \sum_{s \in \mathcal{J}(j',I')} |f_s(x)| &\leq C_\mu D^{-2\mu} \sum_{s \in \mathcal{J}(j',I')} \left(1 + \frac{\operatorname{dist}(x,I_s)}{|I_s|}\right)^{-20} \inf_{y \in I_s} Mf(y) \\
&\leq C_\mu D^{-2\mu+3} \mathbf{1}_{I'}(x) MMf(x).
\end{aligned}
$$

This follows because for a fixed integer $n$, the intervals $I_s$, $s \in \mathcal{J}(j',I')$ with $2^n \leq |I_s| \leq \tfrac{4}{3} 2^n$ can overlap at most $D^3$ times, by (3.46). Hence

$$\sum_{s \in \mathbf{S}} |c_s(x)| \leq C_\mu D^{-2\mu+6} \left\| \sum_{j=1}^{J} \sum_{I \in \mathcal{I}_j} \mathbf{1}_I(x) \right\|_\infty MMf(x),$$

so that the $L^2$ norm is at most

$$\left\| \sum_{s \in \mathbf{S}} |c_s(x)| \right\|_2 \leq C_\mu D^{-2\mu+9}.$$

But $\mu \geq 400$, so (3.50) follows.                                            $\square$

Georgia Institute of Technology, Atlanta, GA
*E-mail address*: lacey@math.gatech.edu
*Homepage*: http://www.math.gatech.edu/~lacey